\documentclass[12pt,reqno]{amsart}
\usepackage{geometry}                
\usepackage{subfig}
\usepackage{amsmath,amssymb}
\usepackage{graphicx}
\usepackage{amssymb,latexsym, amsmath}
\usepackage{amsfonts,amsthm}
\usepackage{amscd}
\usepackage{mathrsfs}
\usepackage{hyperref}
\usepackage{eucal}
\usepackage{multicol}
\usepackage{epstopdf}
\usepackage{amsgen}
\usepackage{xspace}
\usepackage{verbatim}
\usepackage{stmaryrd}
\usepackage{graphicx}
\usepackage{caption}
\usepackage{enumitem}
\newlist{steps}{enumerate}{1}
\setlist[steps, 1]{label = Step \arabic*:}

\newcommand{\gd}{\Delta}

\newcommand{\inpt}[1]{\langle #1 \rangle}

\newcommand{\gw}{\Omega}
\newcommand{\ap}{\alpha}
\newcommand{\ga}{\gamma}

\newcommand{\gb}{\beta}

\newcommand{\gl}{\lambda}

\newcommand{\gs}{\sigma}

\newcommand{\ms}{\mathscr}

\newcommand{\gk}{\kappa}

\newcommand{\nb}{\nabla}
\newcommand{\vp}{\varphi}
\newcommand{\ve}{\varepsilon}
\newcommand{\pdr}{\partial}

\newcommand{\beq}{\begin{equation}}
\newcommand{\eeq}{\end{equation}}
\newcommand{\bea}{\begin{align}}
\newcommand{\eea}{\end{align}}
\newcommand{\bthm}{\begin{theorem}}
\newcommand{\ethm}{\end{theorem}}
\newcommand{\bpr}{\begin{proof}}
\newcommand{\epr}{\end{proof}}
\newcommand{\bcl}{\begin{corollary}}
\newcommand{\ecl}{\end{corollary}}
\newcommand{\bpn}{\begin{proposition}}
\newcommand{\epn}{\end{proposition}}
\newcommand{\bre}{\begin{remark}}
\newcommand{\ere}{\end{remark}}
\newcommand{\bdf}{\begin{definition}}
\newcommand{\edf}{\end{definition}}
\newcommand{\bss}{\begin{align*}}
\newcommand{\ess}{\end{align*}}

\newcommand{\bl}{\label}

\newcommand{\mR}{\mathbb{R}}
\newtheorem{theorem}{Theorem}[section]
\newtheorem{corollary}[theorem]{Corollary}

\newtheorem{proposition}[theorem]{Proposition}

\theoremstyle{definition}
\newtheorem{definition}[theorem]{Definition}
\theoremstyle{remark}
\newtheorem{remark}{Remark}

\numberwithin{equation}{section}

\begin{document}

\title[Memristive Fitz0Hugh-Nagumo Neural Networks]{Synchronization of Memristive FitzHugh-Nagumo Neural Networks}

\author[Y. You]{Yuncheng You}
\address{University of South Florida, Tampa, FL 33620, USA}
\email{you@mail.usf.edu}
\thanks{}

\author[J. Tian]{Jing Tian$^{\dag}$}
\address{Towson University, Towson, MD 21252, USA}
\address{$\dagger$ corresponding author}
\email{jtian@towson.edu$^{\dag}$}
\thanks{}

\author[J. Tu]{Junyi Tu}
\address{Salisbury University, Salisbury, MD 21801, USA}
\email{jxtu@salisbury.edu}
\thanks{}


\date{\today}



\begin{abstract} 
A new mathematical model of neural networks described by diffusive FitzHugh-Nagumo equations with memristors and linear synaptic coupling is proposed and investigated. The existence of absorbing set for the solution semiflow in the energy space is proved and global dynamics of the memristive neural networks are dissipative. Through uniform estimates and maneuver of integral inequalities on the interneuron difference equations, it is shown that exponential synchronization of the neural network at a uniform convergence rate occurs if the coupling strength satisfies a threshold condition explicitly expressed by the system parameters, which is illustrated by an example and numerical simulation experiments. 
\end{abstract}

\maketitle
 
\section{Introduction}

Recently the global dynamics and exponential synchronization of the neural networks modeled by the diffusive Hindmarsh-Rose equations with memristors were proposed and studied by the first author in \cite{Y1, Y2}. In this paper, we shall consider a new mathematical model of neural networks described by the diffusive FitzHugh-Nagumo equations with memristors and linear synaptic interneuron coupling. 

Let a network of $m$ fully coupled memristive neuron cells be denoted by $\mathcal{NW} = \{\mathcal{N}_i : i = 1, 2, \cdots, m\}$, where $m \geq 2$ is a positive integer, which is described by the following model of memristive and diffusive FitzHugh-Nagumo (FHN) equations. Each neuron $\mathcal{N}_i, 1 \leq i \leq m$, in this network is presented by the differential equations:
\beq \bl{cHR}
\begin{split}
	\frac{\pdr u_i}{\pdr t} & = \eta \gd u_i + f (u_i, x) - \gs w_i + J - k \tanh (\rho_i) u_i + \sum_{j = 1}^m P (u_j - u_i),  \\
	\frac{\pdr w_i}{\pdr t} & =  a u_i + c - b w_i,    \\[2pt]
	\frac{\pdr \rho_i}{\pdr t} & = q u_i - r \rho_i,  
\end{split} 
\eeq
for $t > 0,\, x \in \gw \subset \mathbb{R}^{n}$ ($n \leq 3$), where $\gw$ is a bounded domain with locally Lipschitz continuous boundary $\partial \gw$. 

In the membrane potential $u_i$-equations, the nonlinear term $k \tanh (\rho_i) u_i$ presents the memristive coupling effect \cite{RM, W, XQM}, where $\rho_i (t, x)$ stands for the memductance of the memristor and $\tanh (\rho_i)$ is the the electromagnetic induction flux with its coupling strength coefficient $k$. In this system, the fast excitatory variable $u_i(t,x)$ refers to the transmembrane electrical potential of a neuron cell and the slow recovering variable $w_i(t, x)$ represents the integrated ionic current across the neuron membrane. The network neuron coupling terms are assumed to be linear with a common strength coefficient $P$ in the membrane potential equation.

We impose the homogeneous Neumann boundary condition
\begin{equation} \label{nbc}
	\frac{\pdr u_i}{\pdr \nu} (t, x) = 0, \quad \text{for} \;\; t > 0,  \; x \in \partial \gw, \quad 1 \leq i \leq m.
\end{equation}
The initial states of the system \eqref{cHR} will be denoted by 
\begin{equation} \bl{inc}
	 u_i^0 (x) = u_i(0, x), \; w_i^0 (x) = w_i (0, x), \; \rho_i^0 = \rho_i (0, x), \;\; 1 \leq i \leq m.
\end{equation}
The following Assumption is made on the scalar function $f \in C^1 (\mathbb{R} \times \gw)$:
\beq \bl{Asp}
	\begin{split}
	&f(s, x) s \leq  - \gl |s|^4 + \vp (x), \quad s \in \mathbb{R}, \; x \in \gw, \\
	&\frac{\pdr f}{\pdr s} (s, x) \leq \gb,  \quad s \in \mathbb{R}, \; x \in \gw, 
	\end{split}
\eeq
where $\gl$ and $\gb$ are positive constants, $\vp \in L^2 (\gw)$ is a given function. Note that the prototype nonlinear homogeneous function $f(s) = s(s - \kappa)(1 -s)$ with $\kappa > 0$ in the original FitzHugh-Nagumo equations \cite{FH} satisfies the properties \eqref{Asp}, in particular,
\begin{gather*}
	f(s) s = -s^4 + (1 + \kappa)s^3 - \kappa s^2 \leq - \frac{1}{4} s^4 + \frac{1}{4} (1 + \kappa)^4,   \\
	f\,^\prime(s) = - 3s^2 + 2(1 + \gk)s - \gk \leq \frac{1}{3} (1 + \gk)^2 - \gk \leq \beta = \frac{1}{3} (1 + \gk)^2.
\end{gather*}	
All the parameters $\eta, \gs, k, a, c, b, q, r$, and $P$ can be any positive constants, while the reference membrane potential $J$ can be any real number constant. 

Typical models of neuron dynamics are four-dimensional Hodgkin-Huxley equations \cite{HH}, two-dimensional FitzHugh-Nagumo equations \cite{FH}, and three-dimensional Hindmarsh-Rose equations \cite{HR}, which originally consist of ordinary differential equations without memristors and characterize the periodic firing-bursting dynamics for neurons and nerve systems. Analysis through Hopf bifurcations and energy or Hamiltonian methods with semi-numerical simulations are the main approaches to show many solution patterns and collective synchronization behavior \cite{ET, Iz, CPY, ZS}.

Global dynamics and synchronization of ensemble neurons and neural networks modeled by partly diffusive Hindmarsh-Rose equations and FitzHugh-Nagumo equations have been studied by the authors' group in recent years \cite{CPY, PSY, LSY, SkY}. These models are hybrid differential equations and reflect the structural feature of neuron cells, which contain short-branch dendrites receiving incoming signals and long-branch axons propagating and transmitting outgoing signals through synapses. 

The concept of memristor was coined by Chua \cite{Chua} to describe the effect of electromagnetic flux on moving electric charges. Physical and generic memristive systems initially reported in \cite{ChuaK, SS} attracted broad scientific interests in the recent decade, especially recognized in biological neuron models and artificial intelligence computing \cite{Ay, E, Jo, Li, SW, US2} as a new type (other than electrical and chemical) synapsis or as an ideal component which has the nonvolatile properties and can process dynamically memorized signal information to exhibit more complex or chaotic dynamics in neural networks. Memristor-based mathematical models now penetrate many fields with applications to image encryption, DNA sequences operation, brain criticality, cell physiology, cybersecurity, and quantum computers, cf. \cite{Li, LH, RK, Sn, SW, W, WP}. 

The researches on memristive FitzHugh-Nagumo and Hindmarsh-Rose neural networks in ordinary differential equations have been expanding in the recent decade, cf. \cite{Ay, EE, K, RM, ZY} and many references therein. Various synchronization results with memristive effect of these models are achieved \cite{Guan, HY, N, RJ, US2, VK, XJ} mainly by the methods of generalized Hamiltonian functions, Lyapunov exponents, and the computational algebra with numerical simulations.

In this work we shall rigorously prove a threshold condition on the neuron coupling strength $P$ to ensure an exponential synchronization of the memristive neural networks \eqref{cHR} through the approach of dissipative dynamical analysis and sharp uniform estimates, which can be extended to study complex or artificial neural networks.

\section{\textbf{Formulation and Preliminaries}}

Define two Hilbert spaces of functions: 
$$
	E = [L^2 (\gw, \mathbb{R}^3)]^m \quad \text{and}  \quad  \Pi = [H^1 (\gw) \times L^2 (\gw, \mathbb{R}^2)]^m
$$ 
where $H^1 (\gw)$ is a Sobolev space. Call $E$ the energy space and $\Pi$ the regular space. The norm and inner-product of $L^2(\gw)$ or $E$ will be denoted by $\| \, \cdot \, \|$ and $\inpt{\,\cdot , \cdot\,}$, respectively. We use $| \, \cdot \, |$ to denote a vector norm or a set measure in $\mR^n$. The initial-boundary value problem \eqref{cHR}-\eqref{inc} can be formulated into an initial value problem of the evolutionary equation:
\begin{equation} \label{pb}
\begin{split}
	&\frac{\partial g}{\partial t} = A g + F(g), \;\;  t > 0, \\
	&g(0) = g^0 \in E.
\end{split}
\end{equation}
The unknown function is a column vector $g(t) = \text{col}\; (g_1 (t), g_2 (t), \cdots, g_m (t))$, where
$$
	g_i (t) = \text{col}\, (u_i(t, \cdot),\, w_i(t, \cdot),\, \rho_i (t, \cdot)), \quad 1 \leq i \leq m,
$$ 
characterizes the dynamics of the neuron $\mathcal{N}_i$. The initial data function in \eqref{pb} is 
$$
	g(0) = g^0 = \text{col}\; (g_1^0, \,g_2^0, \cdots, g_m^0) \quad \text{where} \;\; g_i^0 = \text{col}\,(u_i^0, \, w_i^0, \,\rho_i^0), \; 1 \leq i \leq m. 
$$
The energy norm $\|g(t)\|$ of the solution for the evolutionary equation \eqref{pb} in the space $E$ is given by 
$$
	\|g(t)\|^2 = \sum_{i=1}^m \|g_i (t)\|^2 = \sum_{i=1}^m \left(\|u_i(t)\|^2 + \|w_i(t)\|^2 + \|\rho_i(t)\|^2 \right).
$$
The closed linear operator $A$ in \eqref{pb} is defined by $A= \text{diag} \, (A_1, A_2, \cdots, A_m)$, where
\begin{equation} \label{opA}
A_i =
\begin{pmatrix}
\eta \gd \quad & 0  \quad & 0  \\[3pt]
0 \quad & - b I \quad  & 0 \\[3pt]
0 \quad & 0 \quad & - r I  \\[3pt]
\end{pmatrix}
: \mathcal{D}(A) \rightarrow E, \quad i = 1, 2, \cdots, m,
\end{equation}
with the domain $\mathcal{D}(A) = \{g \in [H^2(\gw) \times L^2 (\gw, \mathbb{R}^2)]^m: \pdr u_i /\pdr \nu = 0,\, 1 \leq i \leq m\}$, is the generator of the $C_0$-semigroup $\{e^{At}\}_{t \geq 0}$ on the space $E$ and $I$ is the identity operator. By the fact that $H^{1}(\gw) \hookrightarrow L^6(\gw)$ is a continuous imbedding for space dimension $n \leq 3$ and by the Assumption \eqref{Asp}, the nonlinear mapping 
\begin{equation} \label{opf}
F(g) =
\begin{pmatrix}
f (u_1, x) - \gs w_1 + J - k \tanh (\rho_1) u_1 + \sum_{j = 1}^m P (u_j - u_1) \\[4pt]
a u_1 + c \\[4pt]
q u_1  \\[4pt]
\vdots\\[4pt]
f (u_m, x) - \gs w_m + J - k \tanh (\rho_m) u_m + \sum_{j = 1}^m P (u_j - u_m) \\[4pt]
a u_m + c \\[4pt]
q u_m  
\end{pmatrix}
: \Pi \longrightarrow E
\end{equation}
is a locally Lipschitz continuous mapping. 

In this work we shall consider the weak solutions \cite[Section XV.3]{CV} of this initial value problem \eqref{pb}.
\begin{definition} \label{D:wksn}
	A $3m$-dimensional vector function $g(t, x)$, where $(t, x) \in [0, \tau] \times \gw$, is called a weak solution to the initial value problem of the evolutionary equation \eqref{pb}, if the following two conditions are satisfied: 
	
	\textup{(i)} $\frac{d}{dt} \langle g, \xi \rangle = \langle Ag, \xi \rangle + \langle F(g), \xi \rangle$ is satisfied for a.e. $t \in [0, \tau]$ and any $\xi \in E^* = E$. 
		
	\textup{(ii)} $g(t, \cdot) \in  C([0, \tau]; E) \cap C^1 ((0, \tau); E)$ and $g(0) = g^0$.	
	
\noindent Here $E^*$ is the dual space of the Hilbert space $E$.
\end{definition}

The following proposition can be proved by the Galerkin approximation method \cite{CV} with the regularity property \cite{SY} of the parabolic operator semigroup $e^{At}$.

\begin{proposition} \label{pps}
	For any initial state $g^0 \in E$, there exists a unique weak solution $g(t; g^0), \, t \in [0, \tau]$, for some $\tau > 0$ may depending on $g^0$, of the initial value problem \eqref{pb} formulated from the memristive FitzHugh-Nagumo equations \eqref{cHR}. The weak solution $g(t; g^0)$ continuously depends on the initial data and satisfies 
\begin{equation} \label{soln}
	g \in C([0, \tau]; E) \cap C^1 ((0, \tau); E) \cap L^2 ((0, \tau); \Pi).
\end{equation}
Moreover, for any initial state $g^0 \in E$, the weak solution $g(t; g^0)$ becomes a strong solution for $t \in (0, \tau)$, which has the regularity
\begin{equation} \bl{ss}
	g \in C((0, \tau]; \Pi) \cap C^1 ((0, \tau); \Pi).
\end{equation}
\end{proposition}

An infinite dimensional dynamical systems \cite{CV, SY} for time $t \geq 0$ only is called a semiflow. Absorbing set defined below is the key concept to characterize dissipative dynamics of a semiflow.

\begin{definition} \label{Dabsb}
	Let $\{S(t)\}_{t \geq 0}$ be a semiflow on a Banach space $\ms{X}$. A bounded set $B^*$ of $\ms{X}$ is called an absorbing set of this semiflow, if for any given bounded set $B \subset \ms{X}$ there exists a finite time $T_B \geq 0$ depending on $B$, such that $S(t)B \subset B^*$ for all $t  > T_B$. The semiflow is called dissipative if there exists an absorbing set.
\end{definition}

The Young's inequality in a generic form will be used throughout this paper: For any two positive numbers $x$ and $y$, if $\frac{1}{p} + \frac{1}{q} = 1$ and $p > 1, q > 1$, one has
\beq \bl{Yg}
	x\,y \leq \frac{1}{p} \ve x^p + \frac{1}{q} C(\ve, p)\, y^q \leq \ve x^p + C(\ve, p)\, y^q, \quad C(\ve, p) = \ve^{-q/p},
\eeq
where constant $\ve > 0$ can be arbitrarily small. In Section 4, the Gagliardo-Nirenberg interpolation inequalities \cite[Theorem B.3]{SY} will be used in a crucial step to prove the main result on exponential synchronization of the memristive neural networks.

\section{\textbf{Dissipative Dynamics of Memristive FitzHugh-Nagumo Semiflow}}

In this section, we shall prove the global existence of weak solutions in time for the initial value problem \eqref{pb} to establish a solution semiflow of the memristive FitzHugh-Nagumo neural network modeled by \eqref{cHR}.  Then we show the dissipative dynamics in terms of the existence of an absorbing set of this semiflow in the state spaces $E$.

\begin{theorem} \label{T1}
	For any initial state $g^0 \in E$, there exists a unique global weak solution in time, $g(t; g^0) = \textup{col}\, (u_i(t), w_i(t), \rho_i(t): 1 \leq i \leq m), \, t \in [0, \infty)$, to the initial value problem \eqref{pb} of the memristive FitzHugh-Nagumo equations \eqref{cHR} for the neural network $\mathcal{NW}$. 
\end{theorem}

\begin{proof}
Conduct the $L^2$ inner-products of the $u_i$-equation with $C_1 u_i(t,x)$ for $1 \leq i \leq m$, with the scaling constant $C_1 > 0$ to be chosen later. Then sum them up to get		
\begin{equation} \bl{ui}
	\begin{split}
	&\frac{C_1}{2} \frac{d}{dt} \sum_{i = 1}^m \|u_i\|^2 + C_1 \eta\, \sum_{i=1}^m \|\nb u_i\|^2  = -  C_1P \, \sum_{i=1}^m \sum_{j=1}^m \int_\gw ( u_i - u_j)^2\, dx \\
	&+ C_1 \sum_{i=1}^m \int_\gw (f(u_i, x) u_i - \gs u_i w_i + Ju_i - k \tanh (\rho_i) u_i^2) \, dx  \\
	\leq &\, C_1 \sum_{i=1}^m \int_\gw \left[ - \gl |u_i |^4 + |\vp (x)| - \gs u_i w_i + Ju_i + k |\tanh (\rho_i)| u_i^2 \right] dx \\
	\leq &\, C_1 \sum_{i=1}^m \int_\gw \left[ - \gl u_i^4 + |\vp (x) | + \frac{1}{2}\left(\gl u_i^2 + \frac{\gs^2}{\gl} w_i ^2 \right) + \frac{1}{2} \left(\frac{J^2}{\gl} + \gl u_i^2\right) + k u_i^2 \right] dx \\
	= &\, C_1 \sum_{i=1}^m \int_\gw ((\gl + k)u_i^2 - \gl u_i^4)\, dx + \frac{C_1 \gs^2}{2\gl}\sum_{i=1}^m \|w_i \|^2  + C_1 m \left(\|\vp \|_{L^1} + \frac{J^2}{2\gl} |\gw| \right) \\
	\leq &\, - \frac{1}{2} C_1 \gl \sum_{i=1}^m \int_\gw u_i^4 (t, x)\, dx + \frac{C_1 \gs^2}{2\gl}\sum_{i=1}^m \|w_i \|^2 \\
	&\, + C_1 m \left(\|\vp \| |\gw|^{1/2} + \frac{1}{2\gl} ((\gl + k)^2 + J^2)|\gw| \right),
	\end{split}
\end{equation}
where the Gauss divergence theorem and the Assumption \eqref{Asp} are used. In the last step above, we notice that $(\gl + k)u_i^2 \leq \frac{(\gl + k)^2}{2\gl} + \frac{\gl}{2} u_i^4$. Then summing up the $L^2$ inner-products of the $w_i$-equations with $w_i(t, x)$ for $1 \leq i \leq m$, by Young's inequality \eqref{Yg}, we get
\beq \bl{wi}
	\begin{split}
	&\frac{1}{2} \frac{d}{dt} \sum_{i=1}^m \| w_i \|^2 = \sum_{i=1}^m \int_\gw (a u_i w_i + c w_i - b w_i^2) \, dx  \\
	\leq &\, \sum_{i=1}^m \int_\gw  \left[\left(\frac{a^2}{b} u_i^2 + \frac{1}{4} b \,w_i^2\right)  + \left(\frac{c^2}{b} + \frac{1}{4}  b \,w_i^2\right) - b \,w_i^2\right] dx \\
	= &\, \sum_{i=1}^m \frac{a^2}{b} \int_\gw u_i^2 (t, x)\, dx - \frac{b}{2} \sum_{i=1}^m \|w_i\|^2 + \frac{mc^2}{b} |\gw |.
	\end{split}
\eeq
Next, we sum up the $L^2$ inner-products of the $\rho_i$-equations with $\rho_i(t, x),1 \leq i \leq m,$ to obtain
\beq \bl{rhoi}
	\frac{1}{2} \frac{d}{dt} \sum_{i=1}^m \| \rho_i \|^2 = \sum_{i=1}^m \int_\gw (q u_i \rho_i - r \rho_i^2)\, dx \leq \frac{q^2}{2r} \sum_{i=1}^m u_i^2(t, x)\,dx - \frac{r}{2} \sum_{i=1}^m \|\rho_i \|^2. 
\eeq

Now add the above three inequalities. We come up with 
\beq \bl{uw}
	\begin{split}
	& \frac{1}{2}\frac{d}{dt} \sum_{i = 1}^m \left(C_1 \|u_i\|^2 + \| w_i \|^2 + \| \rho_i \|^2\right) + C_1 \eta \sum_{i=1}^m \|\nb u_i\|^2 + C_1 P \sum_{i=1}^m \sum_{j=1}^m \int_{\gw} ( u_i - u_j)^2 dx \\
        \leq &\, - \sum_{i=1}^m \int_\gw \left(\frac{1}{2} C_1 \gl u_i^4 - \left(\frac{a^2}{b} + \frac{q^2}{2r} \right) u_i^2 \right) dx + \sum_{i=1}^m \left( \frac{C_1 \gs^2}{2\gl} - \frac{b}{2}\right) \|w_i \|^2 - \frac{r}{2} \sum_{i=1}^m \|\rho_i \|^2 \\
        & + C_1 m \left(\|\vp \| |\gw|^{1/2} + \frac{1}{2\gl} ((\gl + k)^2 + J^2) |\gw| \right) + \frac{m c^2}{b} |\gw |,  \quad t \in I_{max} = [0, T_{max}),
	\end{split}
\eeq
where $I_{max}$ is the maximal existence interval of a weak solution. Now we can choose the scaling constant to be
\beq \bl{C1}
	C_1 = \frac{b\gl}{2\gs^2} \quad \text{so that} \quad \frac{C_1 \gs^2}{2\gl} - \frac{b}{2} = - \frac{b}{4}\, .
\eeq
With this choice, from \eqref{uw} it follows that
\beq \bl{u1w}
	\begin{split}
	& \frac{1}{2}\frac{d}{dt} \sum_{i = 1}^m \left(C_1 \|u_i\|^2 +  \| w_i \|^2  + \|\rho_i \|^2\right) + C_1 \eta \sum_{i=1}^m \|\nb u_i\|^2 + \frac{b}{4} \sum_{i=1}^m \|w_i \|^2 + \frac{r}{2} \sum_{i=1}^m \|\rho_i \|^2 \\
	& + \sum_{i=1}^m \int_\gw \left(\frac{1}{2} C_1 \gl u_i^4 - \left(\frac{a^2}{b} + \frac{q^2}{2r} \right) u_i^2 \right) dx + C_1 P \sum_{i=1}^m \sum_{j=1}^m \int_{\gw} ( u_i - u_j)^2\, dx \\
	\leq &\, C_1 m \|\vp \|^2 + m \left(C_1 + \frac{C_1}{2\gl} ((\gl + k)^2 + J^2) + \frac{c^2}{b}\right) |\gw |,  \quad  t \in I_{max} = [0, T_{max}).
	\end{split}
\eeq
By completing square and \eqref{C1}, we have
\beq \bl{kiq}
	\begin{split}
	&\sum_{i=1}^m \int_\gw \left[\frac{1}{2} C_1 \gl u_i^4 -  \left(\frac{a^2}{b} + \frac{q^2}{2r} \right) u_i^2 \right] dx \\
	= & \sum_{i=1}^m  \left(\frac{a^2}{b} + \frac{q^2}{2r} \right) \|u_i \|^2 + \sum_{i=1}^m \int_\gw \left[\frac{b \gl^2}{4\gs^2} u_i^4 - 2 \left(\frac{a^2}{b} + \frac{q^2}{2r} \right) u_i^2 \right] dx \\
	= &\sum_{i=1}^m  \left(\frac{a^2}{b} + \frac{q^2}{2r} \right) \|u_i \|^2  + \sum_{i=1}^m \int_\gw \left[ \frac{\sqrt{b} \gl}{2\gs} u_i^2 - \frac{2\gs}{\sqrt{b}\gl} \left(\frac{a^2}{b} + \frac{q^2}{2r}\right) \right]^2 dx \\
	- & \frac{4m \gs^2}{b\gl^2} \left[\frac{a^2}{b} + \frac{q^2}{2r}\right]^2 |\gw | \geq \sum_{i=1}^m  \left(\frac{a^2}{b} + \frac{q^2}{2r} \right) \|u_i \|^2  - \frac{4m \gs^2}{b\gl^2} \left[\frac{a^2}{b} + \frac{q^2}{2r}\right]^2 |\gw |.
	\end{split}
\eeq
Substitute \eqref{kiq} in \eqref{u1w}. It yields the inequality 
\beq \bl{u2w}
	\begin{split}
	& \frac{1}{2}\frac{d}{dt} \sum_{i = 1}^m \left(C_1 \|u_i\|^2 +  \| w_i \|^2  + \|\rho_i \|^2\right) + C_1 \eta \sum_{i=1}^m \|\nb u_i\|^2 + C_1 P \sum_{i=1}^m \sum_{j=1}^m \int_{\gw} ( u_i - u_j)^2\, dx \\
	& + \sum_{i=1}^m  \left(\frac{a^2}{b} + \frac{q^2}{2r} \right) \|u_i \|^2 + \frac{b}{4} \sum_{i=1}^m \|w_i \|^2 + \frac{r}{2} \sum_{i=1}^m \|\rho_i \|^2 \\
	\leq &\, C_1 m \|\vp \|^2 + m \left[C_1 + \frac{C_1}{2\gl} ((\gl + k)^2 + J^2) + \frac{c^2}{b} + \frac{4\gs^2}{b\gl^2} \left[\frac{a^2}{b} + \frac{q^2}{2r}\right]^2 \right] |\gw |, \quad  t \in I_{max}.
	\end{split}
\eeq
Denote by
\beq \bl{C2}
	C_2 = C_1 + \frac{C_1}{2\gl} ((\gl + k)^2 + J^2) + \frac{c^2}{b} + \frac{4\gs^2}{b\gl^2} \left[\frac{a^2}{b} + \frac{q^2}{2r}\right]^2.
\eeq
Then \eqref{u2w} gives rise to the Gronwall-type differential inequality 
\beq \bl{GZ}
	\begin{split}
	 &\frac{d}{dt} \sum_{i = 1}^m \left[ C_1\|u_i\|^2 + \| w_i \|^2 + \|\rho_i \|^2 \right] +\mu \sum_{i = 1}^m \left[C_1\|u_i\|^2 + \| w_i \|^2 + \|\rho_i \|^2 \right]  \\
	\leq &\, \frac{d}{dt} \sum_{i = 1}^m \left[ C_1\|u_i\|^2 + \| w_i \|^2 + \|\rho_i \|^2 \right] + 2 \sum_{i=1}^m \left[ \left(\frac{a^2}{b} + \frac{q^2}{2r} \right) \|u_i \|^2 + \frac{b}{4} \|w_i \|^2 + \frac{r}{2} \|\rho_i \|^2 \right]\\[3pt]
	\leq &\, 2C_1 m \|\vp \|^2 + 2C_2 m |\gw |, \quad \text{for} \;\;  t \in I_{max} = [0, T_{max}),
	\end{split}
\eeq
where 
$$
	\mu = \min \, \left\{ \frac{2a^2}{b} + \frac{q^2}{r}, \; \frac{b}{2}, \; r \right\}.
$$
We can solve the differential inequality \eqref{GZ} to obtain the following bounding estimate of all the weak solutions on the maximal existence time interval $I_{max}$,
\beq \label{dse}
	\begin{split}
	& \|g(t, g^0)\|^2 = \sum_{i=1}^m \|g_i (t, g_i^0)\|^2 = \sum_{i=1}^m \left(\|u_i (t)\|^2 + \|w_i (t)\|^2 + \|\rho_i (t)\|^2 \right) \\
	\leq &\, \frac{\max \{C_1, 1\}}{\min \{C_1, 1\}}e^{- \mu \,t} \|g^0 \|^2 +  \frac{2m}{\mu \min \{C_1, 1\}} \left( C_1 \|\vp\|^2 + C_2|\gw |\right), \quad t \in [0, \infty).
	\end{split}
\eeq
Here it is shown that $I_{max} = [0, \infty)$ for every weak solution $g(t, g^0)$ because it will never blow up at any finite time. Therefore, for any initial state $g^0 = (g_1^0, \cdots, g_m^0) \in E$, there exists a unique global weak solution in time $t \in [0, \infty)$ of the initial value problem \eqref{pb} for this memristive neural network model \eqref{cHR} in the space $E$.
\end{proof}

Based on the global existence of weak solutions shown in Theorem \ref{T1}, we define the solution semiflow $\{S(t): E \to E\}_{t \geq 0}$ of the memristive and diffusive FitzHugh-Nagumo equations \eqref{cHR} to be
$$
	S(t): g^0 \longmapsto g(t; g^0) = \text{col}\, (u_i (t, \cdot), w_i (t, \cdot), \rho_i (t, \cdot): 1 \leq i \leq m), \quad t \geq 0.
$$
We call this semiflow $\{S(t)\}_{t \geq 0}$ the \emph{memristive FitzHugh-Nagumo neural network semiflow} generated by the neural network model equations \eqref{cHR}. 

The next theorem shows that the memristive FitzHugh-Nagumo neural network semiflow $\{S(t)\}_{t \geq 0}$ is a dissipative dynamical system in the state space $E$.
\begin{theorem} \label{Eab}
	There exists a bounded absorbing set for the memristive FitzHugh-Nagumo neural network semiflow $\{S(t)\}_{t \geq 0}$ in the state space $E$, which is the bounded ball 
\beq \label{Br}
	B^* = \{ h \in E: \| h \|^2 \leq K\}.
\eeq 
	Here the constant 
\beq \bl{K}
	K = 1 +  \frac{2m}{\mu \min \{C_1, 1\}} \left( C_1 \|\vp\|^2 + C_2|\gw |\right),
\eeq
where the constants $C_1$ and $C_2$ are given in \eqref{C1} and \eqref{C2}.
\end{theorem}

\begin{proof}
This is the consequence of the global uniform estimate \eqref{dse} shown in the proof of Theorem \ref{T1}, which implies that
\beq \label{lsp}
	\limsup_{t \to \infty} \|g(t; g^0)\|^2 = \limsup_{t \to \infty} \, \sum_{i=1}^m \|g_i(t; g_i^0)\|^2 < K 
\eeq
for all weak solutions of \eqref{pb} with any initial data $g^0$ in $E$. Moreover, for any given bounded set $B = \{h \in E: \|h \|^2 \leq L\}$ in $E$, there exists a finite time 
$$
	T_B = \frac{1}{\mu} \log^+ \left(L\,\frac{\max \{C_1, 1\}}{\min \{C_1, 1\}}\right)
$$
such that all the solution trajectories started at the initial time $t = 0$ from the set $B$ will permanently enter the bounded ball $B^*$ shown in \eqref{Br} for $t > T_B$.  Therefore, the bounded ball $B^*$ is an absorbing set in $E$ for the semiflow $\{S(t)\}_{t \geq 0}$ so that this memristive FitzHugh-Nagumo neural network semiflow is dissipative.
\end{proof}

We shall further prove an ultimate uniform bound of the membrane potential functions $\{ u_i (t, x): 1 \leq i \leq m\}$ for all the weak solutions in the higher-order integrable space $L^4(\gw)$, which paves the way for the attempt to achieve the memristive neural network synchronization in the next section.

\begin{theorem} \bl{T4}
	There exists a constant $Q > 0$ such that for any initial data $g^0 \in E$, the $u_i (t)$ components, $1 \leq i \leq m$, of the weak solution $g(t; g^0) = (g_1 (t), \cdots, g_m (t))$ of the initial value problem \eqref{pb} for the memristive FitzHugh-Nagumo neural network $\mathcal{NW}$ satisfies the absorbing property in the space $L^4(\gw)$,
\beq \bl{Lbd}
	\limsup_{t \to \infty}\, \sum_{i=1}^m \|u_i(t)\|^4_{L^4} < 1 + Q. 
\eeq
\end{theorem}

\begin{proof}
Take the $L^2$ inner-product of the $u_i$-equation in \eqref{cHR} with $u_i^3 (t), 1 \leq i \leq m$, and sum them up. By the boundary condition \eqref{nbc} and Assumption \eqref{Asp}, we get
\beq \label{uL4}
	\begin{split}
	&\frac{1}{4}\, \frac{d}{dt} \sum_{i=1}^m \|u_i(t)\|^{4}_{L^{4}} + 3\eta \sum_{i=1}^m \|u_i \nb u_i \|^2_{L^2} \\
	&+ P \sum_{i=1}^m \sum_{j=1}^m \int_{\gw} (u_i - u_j)^2 (u_i^2 + u_i u_j + u_j^2)\, dx \\
	= &\, \sum_{i=1}^m \int_\gw (f(u_i, x) u_i^3 - \gs u_i^3 w_i + J u_i^3 - k \tanh (\rho_i) u_i^4)\, dx \\
	\leq &\, \sum_{i=1}^m \int_\gw (- \gl u_i^6 + u_i^2 \vp (x) - \gs u_i^3 w_i + Ju_i^3 + k u_i^4)\, dx, \quad t > 0.
	\end{split}
\eeq
By Cauchy inequality, it is seen that
\beq \bl{vy} 
	u_i^2 \vp (x) - \gs u_i^3 w_i + Ju_i^3 \leq \frac{1}{6} \gl u_i^4 + \frac{1}{3} \gl u_i^6 + \frac{6}{\gl} \left(\vp^2 (x) + \gs^2 |w_i(t,x)|^2 + J^2 \right).
\eeq
Using Young's inequality \eqref{Yg}, 
\beq \bl{ky}
	k u_i^4 \leq \frac{1}{3} \left(\frac{16\, k^3}{\gl^2}\right) + \frac{2}{3} \left(\frac{\gl}{4} u_i^6\right) \leq \frac{6 k^3}{\gl^2} + \frac{1}{6} \gl u_i^6.
\eeq
Note that $\sum_{i=1}^m \sum_{j=1}^m (u_i - u_j)^2 (u_i^2 + u_i u_j + u_j^2) \geq 0$ always holds and
$$
	u_i^4 \leq \frac{1}{3} + \frac{2}{3} u_i^6 \leq 1 + u_i^6,
$$
so that
\beq \bl{L46}
-\frac{1}{2}{u_i}^6\leq \frac{1}{2}-\frac{1}{2}{u_i}^4,
\eeq 
From \eqref{uL4} wherein we successively use the above inequalities \eqref{vy}, \eqref{ky}, \eqref{L46} and \eqref{lsp}, it follows that 

\beq \label{u3}
	\begin{split}
	&\frac{1}{4} \, \frac{d}{dt} \sum_{i=1}^m \|u_i(t)\|^{4}_{L^{4}} + 3\eta \sum_{i=1}^m \|u_i \nb u_i \|^2 \\
	\leq &\, \sum_{i=1}^m  \left(\gl \int_\gw \left(\frac{1}{6} u_i^4 - \frac{1}{2} u_i^6\right) dx + \frac{6 \gs^2}{\gl} \|w_i (t)\|^2 \right)+ \frac{6 m}{\gl} \left(\|\vp \|^2 + \left(J^2 + \frac{k^3}{\gl}\right) |\gw |\right) \\
	\leq &\, \sum_{i=1}^m  \left( - \frac{\gl}{3} \int_\gw u_i^4 \, dx + \frac{6 \gs^2}{\gl} \|w_i (t)\|^2 \right)+\frac{m\lambda}{2}|\Omega| + \frac{6 m}{\gl} \left( \|\vp \|^2 + \left(J^2 + \frac{k^3}{\gl}\right) |\gw |\right) \\
	\leq &\, - \frac{\gl}{3} \sum_{i=1}^m \|u_i (t) \|^4_{L^4} + \frac{6\gs^2}{\gl} K + m \left[\frac{6}{\lambda} \|\vp \|^2 +\left(\frac{\lambda}{2}+\frac{6}{\lambda}J^2+\frac{6}{\lambda^2}k^3\right)|\Omega|\right], \;\, t > 0.   
	\end{split}
\eeq
Consequently, with the non-negative gradient term removed, \eqref{u3} shows that
\beq \bl{Gu4}
	\begin{split}
        &\frac{d}{dt} \sum_{i=1}^m \|u_i(t)\|^{4}_{L^{4}} + \frac{4\gl}{3} \sum_{i=1}^m \|u_i (t) \|^4_{L^4} \\
        \leq &\,\frac{24\gs^2}{\gl} K + m \left[\frac{24}{\lambda} \|\vp \|^2 +\left(2\lambda+\frac{24}{\lambda}J^2+\frac{24}{\lambda^2}k^3\right)|\Omega|\right], \quad  t > 0.
        \end{split}
\eeq
By the parabolic regularity stated in Proposition \ref{pps}, for any weak solution $g(t; g^0)$ one has $u_i (1) \in H^1 (\gw) \subset L^4 (\gw)$ for $1 \leq i \leq m$. Then the second statement in Proposition \ref{pps} shows that any weak solution has the regularity 
$$
	\sum_{i=1}^m u_i (t) \in C([1, \infty), H^1 (\gw)) \subset C([1, \infty), L^4 (\gw)). 
$$
Apply the Gronwall inequality to \eqref{Gu4}. It results in the bounding estimate of all the $u_i (t)$ components in the space $L^4 (\gw)$ as follows:
\beq \bl{L4B}
        \sum_{i=1}^m \|u_i(t)\|^{4}_{L^4} \leq e^{- \frac{4\gl}{3} (t - 1)} \sum_{i=1}^m \|u_i (1)\|^{4}_{L^4} + Q, \quad \text{for} \;\, t \geq 0,
\eeq
where the constant $Q$ is independent of any initial data and given by
\beq \bl{Q}
	Q = \frac{18\, \gs^2}{\gl^2} K + m \left[\frac{18}{\lambda^2} \|\vp \|^2 +\left(\frac{3}{2}+\frac{18}{\lambda^2}J^2+\frac{18}{\lambda^3}k^3\right)|\Omega|\right].
\eeq 
Therefore, the claim \eqref{Lbd} of this theorem is proved.
\end{proof}

\section{\textbf{Synchronization of Memristive FitzHugh-Nagumo Neural Networks}} 

In this section, we shall prove the main result on the exponential synchronization of the memristive FitzHugh-Nagumo neural networks described by \eqref{cHR} in the state space $E$. This result provides a sufficient quantitative threshold condition for the neuron coupling strength to reach the neural network synchronization.

\begin{definition}
For a model evolutionary equation of a neural network called NW such as \eqref{pb} formulated from the memristive and diffusive FitzHugh-Nagumo equations \eqref{cHR}, we define the asynchronous degree of this neural network in a state space (as a Banach space) $Z$ to be
$$
	deg_s \,(\text{NW})= \sum_{1\, \leq i \,< j\, \leq \,m} \left\{ \sup_{g_i^0, \, g_j^0\,  \in \, Z} \, \left\{\limsup_{t \to \infty} \, \|g_i (t; g^0_i) - g_j (t; g^0_j)\|_Z \right\}\right\}
$$ 
where $g_i (t)$ and $g_j (t)$ are any two solutions of the model equation with the initial states $g_i^0$ and $g_j^0$ for two neurons $\mathcal{N}_i$ and $\mathcal{N}_j$ in the network. The neural network is said to be asymptotically synchronized if 
$$
	deg_s \,(\text{NW}) = 0.
$$
If the asymptotic convergence to zero of the difference norm above for any two neurons in the network admits a uniform exponential rate, then the neural network is called exponentially synchronized. 
\end{definition}

Introduce the neuron difference functions: For $i, j = 1, \cdots, m$, we define 
$$
	U_{ij} (t,x) = u_i(t,x) - u_j (t,x), \, W_{ij} (t,x) = w_i(t,x) - w_j (t,x), \, R_{ij} (t,x) = \rho_i(t,x) - \rho_j (t,x).
$$
Given any initial state $g^0 = \text{col}\, (g_1^0, \cdots, g_m^0)$ in the space $E$, the difference between any two solutions of \eqref{pb} associated with two neurons $\mathcal{N}_i$ and $\mathcal{N}_j$ in the network is what we consider: 
$$
	g_i (t, g_i^0) - g_j (t, g_j^0) = \text{col}\, (U_{ij}(t, \cdot ), W_{ij}(t, \cdot ), R_{ij}(t, \cdot)), \quad t \geq 0.
$$
By subtraction of the three governing equations for the $j$-th neuron from the corresponding equations for the $i$-th neuron in \eqref{cHR}, we obtain the following differencing FitzHugh-Nagumo equations. For $i, j = 1, \cdots, m$,
\beq \bl{dHR} 
	\begin{split}
	\frac{\pdr U}{\pdr t} = \eta \gd U + f(u_i, x) - & f(u_j, x) - \gs W - k (\tanh (\rho_i)u_i - \tanh (\rho_j) u_j) - mPU,  \\
	&\frac{\pdr W}{\pdr t}  = a U - b W, \\
	&\frac{\pdr R}{\pdr t}  = q U - r R.
	\end{split}
\eeq
Here we can simply write $U(t, x) = U_{ij}(t, x), W (t, x) = W_{ij}(t, x), R(t, x) = R_{ij}(t, x)$ and further $U(t) = U(t, \cdot), W(t) = W(t, \cdot), R(t) = R(t, \cdot)$ for notational convenience. 

The following exponential synchronization theorem is the main result of this paper.
\begin{theorem} \bl{ThM}
	For the memristive FitzHugh-Nagumo neural network $\mathcal{NW}$ with the model \eqref{cHR}-\eqref{Asp}, If the following threshold condition is satisfied by the coupling strength coefficient $P$,  
\beq \bl{SC}
	P > \Gamma.
\eeq
Here the threshold $\Gamma$ is defined to be
$$\Gamma=\frac{1}{m} \left(\gb + k + \frac{1}{2b} |a - \gs|^2 + \frac{q^2}{r} + \frac{C^{*4}\, k^8 (1 + Q)^2}{\eta^3 \,r^4} \right) ,$$
where the constant $Q$ is given in \eqref{Q} and the constant $C^*$ is the coefficient in the Gagliardo-Nirenberg inequality \eqref{intp}, then the neural network $\mathcal{NW}$ is exponentially synchronized in the state space $E$ at a uniform exponential rate $\ap (P)$: 
\beq \bl{rate}
	\ap (P) = \min \left\{b, \, r, \, \left[ 2mP - 2 \left(\gb + k + \frac{1}{2b} |a - \gs|^2 + \frac{q^2}{r} + \frac{C^{*4}\, k^8(1 +Q)^2}{\eta^3 \,r^4} \right) \right] \right\}.
\eeq
\end{theorem}

\begin{proof}
The proof will go through two steps. 

Step 1. Take the $L^2$ inner-products of the first equation in \eqref{dHR} with $U(t)$, the second equation in \eqref{dHR} with $W(t)$, and the third equation in \eqref{dHR} with $R(t)$. Then sum them altogether and use the Assumption \eqref{Asp} to get
\beq \bl{eG} 
	\begin{split}
	&\frac{1}{2} \frac{d}{dt} (\|U (t)\|^2 + \|W (t)\|^2 + \|R (t)\|^2) + \eta \|\nb U(t)\|^2 + b\,\|W(t)\|^2 + r \|R(t)\|^2  \\
	= &\, \int_\gw (f(u_i, x) - f(u_j, x)) U\, dx + \int_\gw [(a - \gs)UW + qUR]\, dx \\
	&\, -\int_\gw k (\tanh (\rho_i)u_i - \tanh (\rho_j) u_j) U]\, dx - mP\|U(t)\|^2 \\
	\leq &\, \int_\gw  \frac{\pdr f}{\pdr s} \left(\ell u_i + (1- \ell) u_j, x \right) U^2\, dx + \int_\gw [(a - \gs)UW + qUR]\,dx \\
	&\, - k \int_\gw \left[\text{sech}^2 (\xi \rho_i + (1 - \xi)\rho_j) R\, u_i U + \tanh (\rho_j)U^2 \right] dx - mP\|U(t)\|^2 \\
	\leq &\,\gb \|U\|^2 + \int_\gw [(a - \gs)UW + qUR]\,dx - k \int_\gw R u_j U\, dx + (k - mP) \|U\|^2, 
	\end{split}
\eeq
where the properties $|\tanh (\rho_j)| \leq 1$ and $\text{sech}^2 (\xi \rho_i + (1 - \xi)\rho_j) \leq 1$ for the hyperbolic functions and $\xi, \ell \in [0,1]$ are used. 

Next we have to treat the two integral terms on the right-hand side of the differential inequality \eqref{eG}. By the Young's inequality \eqref{Yg}, we have

\beq \bl{WR}
	\begin{split}
	&\int_\gw [(a - \gs)UW + qUR]\,dx \\
	\leq &\, \int_\gw \left[\frac{b}{2} W^2(t, x) + \frac{1}{2b}|a - \gs|^2 U^2(t, x)\right] dx + \int_\gw \left[\frac{r}{4} R^2(t, x) + \frac{q^2}{r} U^2(t, x) \right] dx \\
	= &\, \frac{b}{2} \|W(t)\|^2 + \frac{r}{4} \|R(t)\|^2 + \left(\frac{1}{2b} |a - \gs|^2 + \frac{q^2}{r}\right) \|U(t)\|^2, \quad t > 0.
	\end{split}
\eeq
For the last integral term in \eqref{eG}, by the H\"{o}lder inequality, we get
\beq \bl{RU}
	\begin{split}
	- k \int_\gw R u_j U\, dx \leq &\, k \int_\gw \left(\frac{r}{4k} R^2(t, x) + \frac{k}{r} u_j^2 (t, x) U^2(t, x)\right) dx \\
	\leq &\, \frac{r}{4} \|R(t)\|^2 +  \frac{k^2}{r} \left[\int_\gw u_j^4 (t, x)\,dx\right]^{1/2} \left[\int_\gw U^4 (t, x)\,dx\right]^{1/2} \\
	= &\, \frac{r}{4} \|R(t)\|^2 +  \frac{k^2}{r} \|u_j(t)\|^2_{L^4} \|U(t)\|^2_{L^4}, \quad t > 0.
	\end{split}
\eeq
Substitute the term estimates \eqref{WR} and \eqref{RU} into the differential inequality \eqref{eG}, we obtain
\begin{equation} \bl{mG}
	\begin{split}
	&\frac{1}{2} \frac{d}{dt} (\|U (t)\|^2 + \|W (t)\|^2 + \|R (t)\|^2) + \eta \|\nb U(t)\|^2 + \frac{b}{2} \|W(t)\|^2 + \frac{r}{2} \|R(t)\|^2   \\
	\leq &\,\left(\gb + k + \frac{1}{2b} |a - \gs|^2 + \frac{q^2}{r}- mP\right) \|U(t)\|^2 + \frac{k^2}{r} \|u_j(t)\|^2_{L^4} \|U(t)\|^2_{L^4}, \;\,  t > 0.
	\end{split}
\end{equation}

Step 2. The key challenge here is to handle the last term on the right-hand side of the estimate inequality \eqref{mG}. We shall exploit the sharp technique of Gagliardo-Nirenberg interpolation inequalities \cite[Theorem B.3]{SY}. In view of the Sobolev embedding
$$
	H^1 (\gw) \subset L^4 (\gw) \subset L^2 (\gw),
$$
one has
\beq \bl{intp}
	\|U(t) \|^2_{L^4} \leq C^* \|\nb U(t)\|^{2\theta} \|U(t)\|^{2(1 - \theta)}
\eeq
where the coefficient $C^*(\gw) > 0$ only depends on the spatial domain $\gw$, and the interpolation index $\theta = 3/4$ is determined by 
$$
	- \frac{3}{4} = \theta \left(1 - \frac{3}{2} \right) - (1 - \theta)\, \frac{3}{2} \, .
$$
Hence \eqref{intp} shows that 
\beq \bl{intU}
	\|U(t) \|^2_{L^4} \leq C^* \|\nb U(t)\|^{3/2} \|U(t)\|^{1/2}.
\eeq
According to Theorem \ref{T4} and \eqref{Lbd}, $\limsup_{t \to \infty} \sum_{i=1}^m \|u_i(t)\|^4_{L^4} < 1 + Q$, so that there exists a finite time $T(g^0) \geq 0$ such that for all $1 \leq i \leq m$,
$$
	\|u_i (t)\|^2_{L^4} < (1 + Q)^{1/2}, \quad \text{for all} \;\; t > T(g^0). 
$$
Therefore, for any given initial state $g^0 \in E$, by \eqref{intU} and using Young's inequality \eqref{Yg}, we can estimate 
\beq \bl{U42}
	\begin{split}
	&\frac{k^2}{r} \|u_j(t)\|^2_{L^4} \|U(t)\|^2_{L^4} \leq \frac{k^2}{r}(1 + Q)^{1/2} \, \|U(t)\|^2_{L^4} \\
	\leq &\, C^* \|\nb U(t)\|^{3/2} \left[\frac{k^2}{r} (1 + Q)^{1/2} \, \|U(t)\|^{1/2} \right] \\
	\leq &\, \eta \|\nb U(t)\|^{(3/2) \times (4/3)} + \frac{1}{\eta^3} \left[\frac{C^* k^2}{r} (1 + Q)^{1/2} \, \|U(t)\|^{1/2} \right]^4 \\
	= &\, \eta \|\nb U(t)\|^2 + \frac{C^{*4}\, k^8 (1 + Q)^2}{\eta^3 \,r^4} \|U(t)\|^2, \quad t > T(g^0). 
	\end{split}
\eeq
Substitute \eqref{U42} in \eqref{mG}. After cancellation of the gradient terms $\eta \|\nb U(t)\|^2$ on both sides of the inequality, it follows that for $t > T(g^0)$,
\beq \bl{FG}
	\begin{split}
	&\frac{1}{2} \frac{d}{dt} (\|U (t)\|^2 + \|W (t)\|^2 + \|R (t)\|^2)  + \frac{b}{2} \|W(t)\|^2 + \frac{r}{2} \|R(t)\|^2   \\
	+ &\left[ mP - \left(\gb + k + \frac{1}{2b} |a - \gs|^2 + \frac{q^2}{r} + \frac{C^{*4}\,k^8(1 + Q)^2}{\eta^3 \,r^4} \right) \right] \|U(t)\|^2 \leq 0.
	\end{split}
\eeq
Therefore, the following Gronwall-type inequality holds:
\beq \bl{Grw}
	\begin{split}
	& \frac{d}{dt} (\|U (t)\|^2 + \|W (t)\|^2 + \|R (t)\|^2) + \ap(P) (\|U (t)\|^2 + \|W (t)\|^2 + \|R (t)\|^2)   \\
	\leq &\, \frac{d}{dt} (\|U (t)\|^2 + \|W (t)\|^2 + \|R (t)\|^2) + (b \|W(t)\|^2 + r \|R(t)\|^2)   \\
	+ &\, 2\left[ mP - \left(\gb + k + \frac{1}{2b} |a - \gs|^2 + \frac{q^2}{r} + \frac{C^{*4}\, k^8(1 + Q)^2}{\eta^3 \,r^4} \right) \right] \|U(t)\|^2 \leq 0, 
	\end{split}
\eeq
for $t > T(g^0)$, where the constant 
$$
	\ap (P) = \min \left\{b, \, r, \, \left[ 2mP - 2 \left(\gb + k + \frac{1}{2b} |a - \gs|^2 + \frac{q^2}{r} + \frac{C^{*4}\, k^8 (1 + Q)^2}{\eta^3 \,r^4} \right) \right] \right\}
$$
as shown in \eqref{rate}.

Under the threshold condition \eqref{SC} of this theorem, solving this linear Gronwall inequality \eqref{Grw} directly shows the exponential synchronization result: For any initial state $g^0 \in E$ and any two neurons $\mathcal{N}_i$ and $\mathcal{N}_j$ in the memristive FitzHugh-Nagumo neural network \eqref{cHR}, their difference function $g_i(t; g_i^0) - g_j(t; g_j^0)$ converges to zero in the state space $E$ exponentially at a uniform rate $\ap (P)$ shown in \eqref{rate}. Namely, for any $1 \leq i <  j \leq m$,
\beq \bl{ESyn} 
	\begin{split}
	\| g_i (t) - g_j (t) \|_E^2 &= \|U_{ij}(t)\|^2 + \|W_{ij}(t)\|^2 + \|R_{ij}(t)\|^2  \\[2pt]
	&\leq e^{- \ap (P) \,t } \left\|g_i^0 - g_j^0 \right\|^2  \to 0, \;\; \text{as} \;\, t \to \infty. 
	\end{split}
\eeq
Hence it is proved that
\beq \bl{degs}
	deg_s (\mathcal{NW}) = \sum_{1 \,\leq \,i  \,<  \,j \,\leq \,m} \left\{\sup_{g^0\, \in \, E} \, \left\{\limsup_{t \to \infty} \|g_i (t) -g_j(t) \|^2_E \right\}\right\} = 0.
\eeq
Thus the exponential synchronization of the memristive and diffusive Hindmarsh-Rose neural network $\mathcal{NW}$ in the space $E$ is proved. 
\end{proof}

\section{\textbf{Example and Numerical Simulation}}

In this section, we test some numerical experiments to verify and illustrate the obtained theoretical result on synchronization stated in Theorem \ref{ThM}.

We numerically solve the memristive FitzHugh-Nagumo neural network $\mathcal{NW}$ with the model \eqref{cHR}-\eqref{Asp} in a two-dimensional square domain. We use the finite difference method for the numerical scheme and programmed in Python.

\begin{figure}%
	\centering
	\subfloat{{\includegraphics[width=8cm]{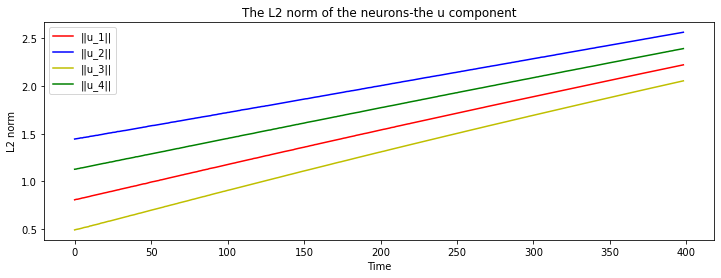} }}%
	\qquad
	\subfloat{{\includegraphics[width=8cm]{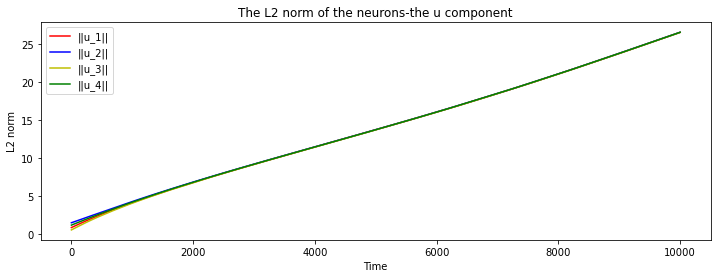} }}%
	\caption{The $L^2$ norm of the neurons $u_i$ at the beginning (the upper figure) and after 10000 iterations (the lower figure)}%
	\label{fig1}%
\end{figure}

\begin{figure}%
	\centering
	\subfloat{{\includegraphics[width=8cm]{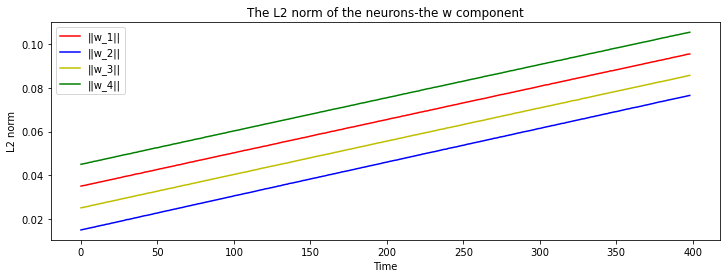} }}%
	\qquad
	\subfloat{{\includegraphics[width=8cm]{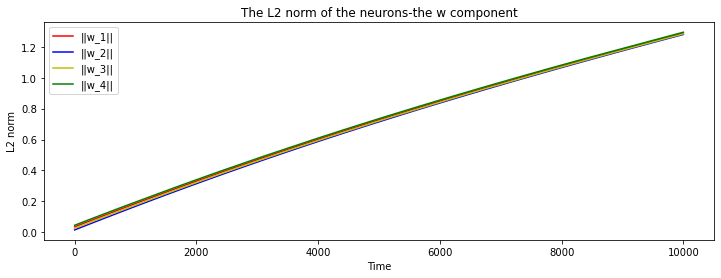} }}%
	\caption{The $L^2$ norm of the neurons $w_i$ at the beginning (the upper figure) and after 10000 iterations (the lower figure)}%
	\label{fig2}%
\end{figure}

\begin{figure}%
	\centering
	\subfloat{{\includegraphics[width=8cm]{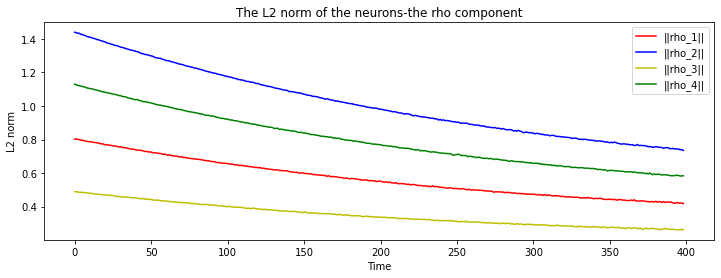} }}%
	\qquad
	\subfloat{{\includegraphics[width=8cm]{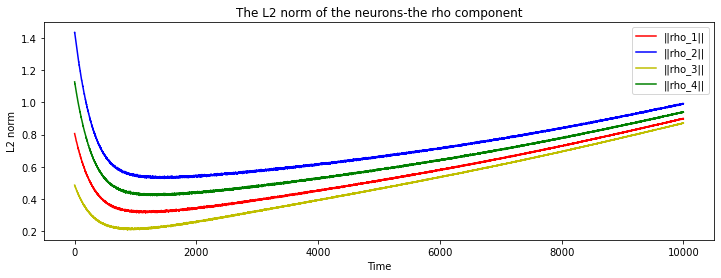} }}%
	\caption{The $L^2$ norm of the neurons $\rho_i$ at the beginning (the upper figure) and after 10000 iterations (the lower figure)}%
	\label{fig3}%
\end{figure}

\begin{figure}%
	\centering
	\subfloat{{\includegraphics[width=8cm]{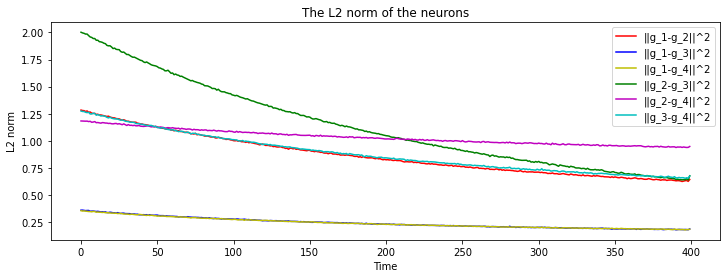} }}%
	\qquad
	\subfloat{{\includegraphics[width=8cm]{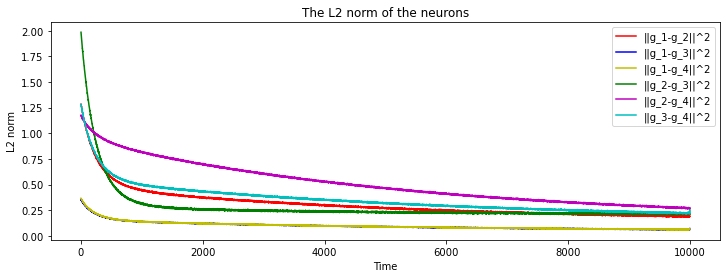} }}%
	\caption{The $L^2$ norm of the neurons from the Energy space at the beginning (the upper figure) and after 10000 iterations (the lower figure)}%
	\label{fig4}%
\end{figure}

Choosing $f(s)=s(s-1)(1-s)$, we consider the following selection of parameters:
\begin{gather*}
	m = 4,\ \eta=10;\ \sigma=0.01;\ J=0.5;\ P=1.45;  \\
	a=0.35;\ b=0.35;\ c=0.7;\ q=0.35;\ r=10.
\end{gather*}
Make the time-step to be 0.00025s and spatial-step to be 1 on a $32*32$ membrane. We compute the $L^2$ norm of neuron potential $u_i$, the recovering variable $w_i$, the memductance $\rho_i$, and also the vector solutions $g_i$ from (\ref{ESyn}) in the energy space $E$ as showing in Figure \ref{fig1} to Figure \ref{fig4}. 

In Figure \ref{fig1} to Figure \ref{fig3}, with a comparison between the beginning stages and results after $10000$ iterations, one can observe the synchronization tendency of the three characterizing variables $(u_i, w_i, \rho_i)$ among the neurons in the simulated mimristive neural network. From Figure \ref{fig4}, we observe that the differences among $\|g_i\|$ tend to $0$. 

We can get the the following constants that used in Theorem \ref{ThM} based on our selection of parameters.
\begin{gather*}
	\lambda =0.25, \quad \phi (x) =4, \quad \beta =\frac{4}{3},  \quad  \\
	C_1 = 437.5\quad C_2 =876.4 \quad \mu =0.175 \quad K =41345645.6 \quad 1 + Q =1220899.6 \quad  C^* = 0.4\\
 	P=1.45 > \Gamma=0.45, \quad \alpha=0.35.
\end{gather*}
Remark: The $\lambda,\ \phi,\ \beta$ from (\ref{Asp}) are not unique. The constant $C^*$ from Gargliardo-Nirenberg inequality is chosen to be $0.4$ based on \cite{Rb}.

Table \ref{tab:table1} to Table \ref{tab:table3} list the sampled values of the three components $u_i, w_i$, and $\rho_i$ of the simulated solution $g_i$ at one same point in the domain at $t=0$ and at the $10000$ time-step. It is seen that with a big difference on the initial values, after a certain time, the curves of $u_i$, $w_i$, and $\rho_i$ tend to be close to each other between various neurons. 

\begin{table}[h!]
	\begin{center}
		\caption{Comparision of the $u_i$ at the point $x=10,\ y=10$}
		\label{tab:table1}
		\begin{tabular}{l|c|r} 
			 & Initial Value & At the $10000$ time step\\
			\hline
			$u_1$ & 0.021435249976028286 & 0.9004893259400936\\
			$u_2$ & 0.04741166022718009 & 0.900530237723381\\
			$u_3$ & 0.01014072281752508 & 0.9004848987658973\\
			$u_4$ & 0.032459331285605755 & 0.9004938438064932\\
		\end{tabular}
	\end{center}
\end{table}

\begin{table}[h!]
	\begin{center}
		\caption{Comparision of the $w_i$ at the point $x=10,\ y=10$}
		\label{tab:table2}
		\begin{tabular}{l|c|r} 
			& Initial Value & At the $10000$ time step\\
			\hline
			$w_1$ & 0.038162808368984426 & 1.4968727659565046\\
			$w_2$ & 0.016109788021287225 & 1.4881936913572227\\
			$w_3$ & 0.028702864538319866 & 1.4927071463421802\\
			$w_4$ & 0.04497009552599465 & 1.4999646128986852\\
		\end{tabular}
	\end{center}
\end{table}

\begin{table}[h!]
	\begin{center}
		\caption{Comparision of the $\rho_i$ at the point $x=10,\ y=10$}
		\label{tab:table3}
		\begin{tabular}{l|c|r} 
			& Initial Value & At the $10000$ time step\\
			\hline
			$\rho_1$ & 0.02516293441173103 & 0.03032291951665976\\
			$\rho_2$ & 0.040098450204586085 & 0.030324374854761176\\
			$\rho_3$ & 0.015327673754565209 & 0.030322776082447638\\
			$\rho_4$ & 0.03404629085391771 & 0.03032306624583209\\
		\end{tabular}
	\end{center}
\end{table}

The synchronization result rigorously proved in this work is illustrated by the example with sample selections of the system parameters and a randomized set of initial data. Our numerical simulation also exhibits that the neuron potentials $u_i$ seem to be synchronized fastest within a limited time, while it takes much longer time to observe the synchronization on the other two variables $w_i$ and $\rho_i$. 

This observation actually enhances the conjecture that adding a nonlinear memristor coupling in the neuron potential equation would accelerate the synchronization for the main variable of neuron membrane potential. On the other hand, it also hints that although the main result Theorem \ref{ThM} confirmed the exponential synchronization has a uniform but may be small convergence rate, each of the three components may have a different synchronization rate, which turns out to be a new open and interesting problem for further research.

\section{\textbf{Conclusions}} We summarize the new contribution of results in this paper. 

1. We propose a new mathematical model in \eqref{cHR} of memristive neural networks in terms of the partly diffusive FitzHugh-Nagumo equations with a nonlinear memristor and linear synaptic coupling in the membrane potential equations. This model as a hybrid system of partial-ordinary differential equations features the full synaptic coupling among all the neurons. In comparison with the extensively studied ODE models, this model is more meaningful to capture the structure of biological neuron cells with long-branch axons in neurodynamics. 

2. The nonlinear memristive coupling across the neurons membrane for this work is in the form of $k \tanh (\rho_i) u_i$, which appeared in quite a few researches as cited in the references. In this paper we proved the exponential synchronization of this model of memristive neural networks. There is another type of coupled memristors in the quadratic form $k(c + \ga \rho_i + \delta \rho_i^2)u_i$, which have also been actively studied most with the ODE models of Hindmarsh-Rose equations on various topics of neuromorphic patterns and chaotic dynamics. But synchronization of the hybrid PDE neural network models such as the diffusive FitzHugh-Nagumo equations with quadratic memristors and linear membrane potential coupling seems still an open problem.

3. In this work we take the analytic approach of global dynamics for the weak solutions to pursue the synchronization investigation. Through the uniform \emph{a priori} estimates of grouped component solutions with adjustable scaling and maneuvering the integral inequalities, we are able to show the existence of absorbing set in the $L^2$-energy space of the solution semiflow and the existence of asymptotic ultimate bound in the higher-order integrable space of the key component solutions. This signifies our methodology from dissipative global dynamics to synchronization for such a complex neural network system. 

4. The spirit of the entire mathematical proof is to tackle and control the nonlinear memductance-potential effect by the linear network coupling in different integrable spaces. Many steps of sharp analysis including the crucial Gagliardo-Nirenberg interpolation are carried out and cohesively managed.

5. The main result Theorem \ref{ThM} of this paper provides a sufficient threshold condition for achieving the exponential synchronization of the memristive FitzHugh-Nagumo neural networks described by this hybrid system. Importantly this quantitative threshold condition \eqref{SC}-\eqref{rate} is simply on the linear coupling coefficient $P$ and the exponentially decaying rate is explicitly expressed by the given biological and mathematical parameters. 

It is expected that modeling of biological and artificial neural networks by hybrid differential equations or other types of PDE in neuroscience and in deep learning field can be generalized with new features such as memristors and time delays. The mathematical approach presented in this work can be further explored and extended with applications in a broad scope.

6. Numerical simulation of our model is presented and aligned with the theoretical proof. All the important constants in the proof are calculated and the threshold of coupling strength P is estimated based on the theoretical result. Synchronization of our example can be seen from the $L^2$ norm figures. Visualization of our model as well as other hybrid equation models of network systems and further investigation of time steps and computational estimate of suboptimal thresholds needed to achieve some required synchronization is another interesting research problem.

\textbf{Acknowledgment}
\vspace{3pt}

Jing Tian's work is supported in part by the AMS Simons Travel Grant and the Jess and Mildred Fisher Endowed Professor fund of Mathematics from the Fisher College of Science and Mathematics at Towson University.

\bibliographystyle{amsplain}

\end{document}